\documentclass[review]{elsarticle}          

\usepackage{url}

\usepackage{amssymb,amsmath,amscd}

\catcode`@=11
\def\plslash{\ifx\@currsize\normalsize
{\mathchoice
{\mbox{\raisebox{0.2ex}{$\scriptstyle\circ$}\kern-1ex$\setminus$}}
{\mbox{\raisebox{0.2ex}{$\scriptstyle\circ$}\kern-1ex$\setminus$}}%
{\mbox{\raisebox{0.14ex}{$\scriptscriptstyle\circ$}\kern-0.8ex%
${\scriptstyle\setminus}$}}%
{\,\mbox{\raisebox{0.14ex}{$\scriptscriptstyle\circ$}\kern-0.8ex%
${\scriptstyle\setminus}$}}}%
\else\ifx\@currsize\large\,\mbox{\raisebox{0.2ex}{$\scriptstyle\circ$}\kern-1ex$\setminus$}
\else\ifx\@currsize\small\,\mbox{\raisebox{0.2ex}{$\scriptstyle\circ$}\kern-1ex$\setminus$}
\else\,\mbox{\raisebox{0.2ex}{$\scriptstyle\circ$}\kern-1ex$\setminus$}
\fi\fi\fi}

\def\prslash{\ifx\@currsize\normalsize
{\mathchoice
{\mbox{\raisebox{0.2ex}{$\scriptstyle\circ$}\kern-1ex$/$}}
{\mbox{\raisebox{0.2ex}{$\scriptstyle\circ$}\kern-1ex$/$}}%
{\mbox{\raisebox{0.14ex}{$\scriptscriptstyle\circ$}\kern-0.8ex%
${\scriptstyle/}$}}%
{\mbox{\raisebox{0.14ex}{$\scriptscriptstyle\circ$}\kern-0.8ex%
${\scriptstyle/}$}}}%
\else\ifx\@currsize\large\mbox{\raisebox{0.2ex}{$\scriptstyle\circ$}\kern-1ex$/$}
\else\ifx\@currsize\small\mbox{\raisebox{0.2ex}{$\scriptstyle\circ$}\kern-1ex$/$}
\else\mbox{\raisebox{0.2ex}{$\scriptstyle\circ$}\kern-1ex$/$}
\fi\fi\fi}

\newcommand{\lslash}{\protect\plslash}
\newcommand{\rslash}{\protect\prslash}
\def\pomoins{\ifx\@currsize\normalsize
\mbox{ $\circ\kern-1.48ex-$ }
\else\ifx\@currsize\large\mbox{ $\circ\kern-1.45ex-$ }
\else\ifx\@currsize\small\mbox{ $\circ\kern-1.51ex-$ }
\else\mbox{ $\circ\kern-1.40ex-$ }\fi\fi\fi}

\catcode`@=12

\newcommand{\rfrac}[2]{\,\rule[0.2ex]{0.1ex}{0.4ex}\kern-0.27ex%
\frac{#1\,}{\,#2}\kern-0.29ex\rule[0.55ex]{0.1ex}{0.4ex}\,\mbox{}}

\newcommand{\lfrac}[2]{\,\rule[0.55ex]{0.1ex}{0.4ex}\kern-0.29ex%
\frac{\,#1}{#2\,}\kern-0.27ex\rule[0.2ex]{0.1ex}{0.4ex}\,\mbox{}}

\newcommand{\id}{\mathsf{Id}}
\newcommand{\im}{\mathsf{Im}}

\newcommand{\cod}{\mathcal{C}_{\mathtt{O}}(\mathcal{S})}
\newcommand{\cd}{\mathcal{C}(\mathcal{S})}
\newcommand{\Zset}{\mathbb{Z}}

\newtheorem{thm}{Theorem}
\newtheorem{propriete}[thm]{Property}
\newtheorem{prop}[thm]{Property}
\newtheorem{cor}[thm]{Corrollary}
\newtheorem{notation}[thm]{Notation}
\newtheorem{defn}[thm]{Definition}
\newtheorem{exmp}[thm]{Example}
\newproof{pf}{Proof}

\newtheorem{rem}[thm]{Remark}

\begin{document}

\begin{frontmatter}

\title{Interval systems over idempotent semiring\\
\small{Extended Version of "Interval systems over idempotent semirings ", Linear Algebra and its Applications, vol. 431, n° 5-7, p. 855-862, 2009.}} 

\author{Laurent Hardouin\corref{cor1}}\ead{laurent.hardouin@univ-angers.fr}
\author{Bertrand Cottenceau} \ead{bertrand.cottenceau@univ-angers.fr}
\author{Mehdi Lhommeau}  \ead{mehdi.lhommeau@univ-angers.fr}  
\author{Euriell Le Corronc} \ead{euriell.lecorronc@univ-angers.fr}

\address{Laboratoire d'Ing\'enierie des Syst\`emes Automatis\'es, \\Universit\'e d'Angers,\\ 62, avenue Notre Dame du Lac, 49000 Angers, France.}  
\cortext[cor1]{Corresponding author}


\begin{abstract}                          
This paper deals with solution of
inequality $\textbf{A}\otimes \textbf{x}\preceq \textbf{b}$, where
$\textbf{A}, \textbf{x}$ and $\textbf{b}$ are interval matrices with entries defined over idempotent semiring.
It deals also with the computation of a pair of
intervals, ($\textbf{x},\textbf{y}$) which satisfies the equation
$\textbf{A} \otimes \textbf{x}=\textbf{B}\otimes \textbf{y}$. It
will be shown that this equation may be solved by considering the
interval version of the iterative scheme proposed in
\cite{cuning03}.

\end{abstract}

\end{frontmatter}

\section{Introduction}

Many problems in the optimization theory and other fields of
mathematics are non-linear in the traditional sense but appear to
be linear over idempotent semirings (e.g., see
\cite{BICOQ,carre79,cuning,gaubert06}). Idempotency of the
additive law induces that idempotent semirings are (partially)
ordered sets. The Residuation theory
\citep{blyth72,croisot56,derderian67} is a suitable tool to deal
with inverse problems of order preserving mappings. It is usually
used to solve equations defined over
idempotent semiring \citep{BICOQ,cuning,cuning03}, for instance the greatest solution of inequality $Ax\preceq b$ may be computed by means of this theory.

 Interval mathematics was pioneered by R.E. Moore (see \citep{moore79})
as a tool for bounding rounding errors in computer programs.
Since then, interval mathematics has been developed into a general
methodology for investigating numerical uncertainty in numerous
problems and algorithms. In \citep{litvinov01c} the idempotent
version is addressed. The authors show that idempotent interval
mathematics appears to be remarkably simpler than its traditional
analog. For example, in the traditional interval arithmetic,
multiplication of intervals is not distributive with respect to
addition of intervals, while idempotent interval arithmetic keeps
this distributivity. This paper deals first with solution of
inequality $\textbf{A}\otimes \textbf{x}\preceq \textbf{b}$, where
$\textbf{A}, \textbf{x}$ and $\textbf{b}$ are interval matrices (see proposition \ref{res_interval}). When equality is achieved, according to definition given
in \citep{cechlarova02}, the equations system  is said weakly
solvable since at least one of its subsystems is solvable. In a
second step, the paper deals with the computation of a pair of
intervals, ($\textbf{x},\textbf{y}$) which satisfies the equation
$\textbf{A} \otimes \textbf{x}=\textbf{B}\otimes \textbf{y}$. It
will be shown that this equation may be solved by considering the
interval version of the iterative scheme proposed in
\citep{cuning03}.

\section{Preliminaries}
\begin{defn} \label{def:semiring} A \textit{semiring} $\mathcal{S}$ is a set endowed with two
internal operations denoted by $\oplus$ (addition) and $\otimes$
(multiplication), both associative and both having neutral
elements denoted by $\varepsilon$ and $e$ respectively, such that
$\oplus$ is also commutative and idempotent (\textit{i.e.}
$a\oplus a =a$). The $\otimes$ operation is distributive with
respect to $\oplus$, and $\varepsilon$ is absorbing for the
product (\textit{i.e.} $ \forall a,~~\varepsilon \otimes a = a
\otimes \varepsilon = \varepsilon$). When $\otimes$ is
commutative, the semiring is said to be commutative.
\end{defn}

 Semirings can
be endowed with a canonical \textit{order} defined  by : $ a
\succeq b~$ iff $~a=a \oplus b$. Then they become sup-semilattices
and $a \oplus b$ is the least upper bound of $a$ and $b$. A
semiring is \textit{complete} if sums of infinite number of terms
are always defined, and if multiplication distributes over
infinite sums too. In particular, the sum of all elements of a complete
semiring is defined and denoted by $\top$ (for 'top'). A complete
semiring (sup-semilattice) becomes a complete lattice for which
the greatest lower bound of $a$ and $b$ is denoted $a \wedge b$,
$i.e.$,  the least upper bound of the (nonempty) subset of all
elements which are less than $a$ and $b$ (see
\citep[$\S4$]{BICOQ}).

\begin{exmp}[$(max,+)$ algebra]
The set $\overline{\Zset} = \Zset \cup \{ - \infty, + \infty \}$
endowed with the max operator as $\oplus$ and the classical sum
 as $\otimes$ is a complete idempotent semiring of which $\varepsilon = - \infty$,
$e=0$ and $\top = + \infty$ and the greatest lower bound $a \wedge
b=min(a,b)$.
\end{exmp}
\begin{exmp}[$(max,min)$ algebra]
The set $\overline{\Zset} = \Zset \cup \{ - \infty, + \infty \}$
endowed with the max operator as $\oplus$ and the min operator as $\otimes$ is a complete idempotent semiring
of which $\varepsilon = -\infty$, $e=+\infty$ and $\top = +
\infty$ and the greatest lower bound $a \wedge b=min(a,b)$.
\end{exmp}
\begin{defn}[Subsemiring]
\label{def:subsemiring} A subset $\mathcal{C}$ of a semiring is
called a subsemiring of $\mathcal{S}$ if
\begin{enumerate}
\item[$\bullet$] $\varepsilon \in \mathcal{C}$ and $e \in
\mathcal{C}$ ; \item[$\bullet$] $\mathcal{C}$ is closed for
$\oplus$ and $\otimes$, i.e, $\forall a,b \in \mathcal{C}$, $a
\oplus b \in \mathcal{C}$ and $a \otimes b \in \mathcal{C}$.
\end{enumerate}
\end{defn}

\begin{defn}[Principal order ideal]
 Let $\mathcal{S}$ be an idempotent semiring. An order ideal set
is a nonempty subset $\mathcal{X}$ of $\mathcal{S}$ such that $(x
\in \mathcal{X}$ and $y \preceq x) \Rightarrow y \in \mathcal{X}$.
A principal order ideal (generated by $x$) is an order ideal,
denoted $_\downarrow\mathcal{X}_x$, of the form
$_\downarrow\mathcal{X}_x := \{y \in \mathcal{S} | y \preceq x \}$.
\end{defn}

The residuation theory provides, under some assumptions,
\textit{greatest} solutions to inequalities such as $f(x) \preceq
b$ where $f$ is an order preserving mapping ($i.e.$, $a\preceq b
\Rightarrow f(a) \preceq f(b)$) defined over ordered sets.

\begin{defn}[Residual and residuated mapping] \label{def_residuation}
An order preserving mapping $f : \mathcal{D} \rightarrow
\mathcal{E}$, where $\mathcal{D}$ and $\mathcal{E}$ are ordered
sets, is a \textit{residuated mapping} if for all $y \in
\mathcal{E}$, the least upper bound of the subset $\{x | f(x)
\preceq y\}$ exists and belongs to this subset. It is then denoted
by $f^\sharp(y)$. Mapping $f^\sharp$ is called the residual of
$f$. When $f$ is residuated, $f^\sharp$ is the unique order
preserving mapping such that \small
\begin{equation}\label{residuated1}
\begin{array}{lcl}
f \circ f^\sharp \preceq \id_{\mathcal{E}} & \textnormal{~~~and~~~} &
f^\sharp \circ f \succeq \id_{\mathcal{D}},
\end{array}
\end{equation}
\normalsize where $\id$ is the identity mapping respectively on
$\mathcal{D}$ and $\mathcal{E}$.
\end{defn}

\begin{propriete} \label{prop_fsharp}
Let $f : \mathcal{D} \rightarrow \mathcal{E}$ be a residuated
mapping, then
\small
$$
\begin{array}{lcl}
y \in \im f & \Leftrightarrow & f(f^\sharp(y)) = y.
\end{array}
$$
\normalsize
\end{propriete}

\begin{propriete}[{\citep[Th. $4.56$]{BICOQ}}]
\label{prop:rescompo} If $h : \mathcal{D} \rightarrow \mathcal{C}$
and $f : \mathcal{C} \rightarrow \mathcal{B}$ are residuated
mappings, then $f \circ h$ is also residuated and \small
\begin{equation}
\begin{array}{lcl}
(f \circ h)^\sharp & = & h^\sharp \circ f^\sharp.
\end{array}
\end{equation}
\normalsize
\end{propriete}

\begin{thm}[{\citep[$\S 4.4.2$]{BICOQ}}]
\label{thm:residuation} Consider the mapping $f:\mathcal{E}
\rightarrow \mathcal{F}$ where $\mathcal{E}$ and $\mathcal{F}$ are
complete semirings. Their bottom elements are, respectively,
denoted by $\varepsilon_{\mathcal{E}}$ and
$\varepsilon_{\mathcal{F}}$. Then, $f$ is residuated iff
$f(\varepsilon_{\mathcal{E}}) = \varepsilon_{\mathcal{F}}$ and
$f(\bigoplus_{x \in \mathcal{G}}x) = \bigoplus_{x \in
\mathcal{G}}f(x)$ for each $\mathcal{G} \subseteq \mathcal{E}$
(i.e., $f$ is \textit{lower-semicontinuous}).
\end{thm}

\begin{cor}
\label{cor:eq_res} The mappings $L_a : x \mapsto a \otimes x$ and
$R_a : x \mapsto x\otimes a$ defined over a complete semiring
$\mathcal{S}$ are both residuated. Their residuals are usually
denoted, respectively, by $L^\sharp_a(x)= a \lslash x$ and
$R^\sharp_a(x) = x \rslash a$ in literature. Therefore, $a\lslash
b$ (resp. $b\rslash a$) is the greatest solution of $a\otimes x
\preceq b$ (resp. $x \otimes a \preceq b$). For matrices, the
practical computation is obtained as follows,
\begin{eqnarray}
\left ( A \lslash B \right )_{ij} & = & \bigwedge \limits_{k=1
\ldots
n} \left ( A_{ki} \lslash B_{kj} \right ), \label{MatrixLeftResiduation}\\
\left ( B \rslash C \right )_{ij} & = & \bigwedge \limits_{k=1
\ldots m} \left ( B_{ik} \rslash C_{jk} \right )
\end{eqnarray}
with $A\in \mathcal{S}^{n \times p}$, $B\in \mathcal{S}^{n \times
m}$ and $C\in \mathcal{S}^{p \times m}$.
\end{cor}

\begin{exmp}
Let $A=\begin{pmatrix}1 & 2 \\ 3& 4 \\ 5 & 6\end{pmatrix}$ and
$B=\begin{pmatrix}8 \\ 9 \\ 10\end{pmatrix}$ be matrices with
entries in $(max,+)$ algebra.

 In $(max,+)$ algebra $a_{ij}\lslash
b_j= b_j-a_{ij}$ then the greatest $x$ such that $A \otimes
x\preceq b$ is given by :
$$x= A \lslash B =\begin{pmatrix}
(1\lslash 8) \wedge (3 \lslash 9) \wedge (5 \lslash 10) \\
(2\lslash 8) \wedge (4 \lslash 9) \wedge (6 \lslash 10)
\end{pmatrix} = \begin{pmatrix}5 \\ 4 \end{pmatrix}$$
\end{exmp}
\begin{exmp}\label{exmp:residuationmaxmin}
Let $A=\begin{pmatrix}1 & 2 \\ 3& 4 \\ 5 & 6\end{pmatrix}$ and
$B=\begin{pmatrix}8 \\ 1 \\ 10\end{pmatrix}$ be matrices with
entries in $(max,min)$ algebra. In $(max,min)$ algebra,
 if $b_j\succeq a_{ij}$ then $a_{ij}\lslash b_j= b_j$ else $a_{ij}\lslash b_j=
 \top$.
$$ A \lslash B =\begin{pmatrix}
(1\lslash 8) \wedge (3 \lslash 1) \wedge (5 \lslash 10) \\
(2\lslash 8) \wedge (4 \lslash 1) \wedge (6 \lslash 10)
\end{pmatrix} = \begin{pmatrix}8 \\ 8 \end{pmatrix}$$
\end{exmp}

The problem of mapping restriction and its connection with the
residuation theory is now addressed. 

\begin{prop}[\citep{blyth72}]\label{prop:canonicalinjection}
 Let $\mathsf{Id}_{|\mathcal{S}_{sub}} : \mathcal{S}_{sub} \rightarrow
 \mathcal{S}$, $x \mapsto x$ be the canonical injection from a complete subsemiring into a
 complete semiring. The injection $\mathsf{Id}_{|\mathcal{S}_{sub}}$
 is residuated and its residual is a projector which will be
 denoted by $\mathsf{Pr}_{\mathcal{S}_{sub}}$, with :
 \small
$$
\begin{array}{lclcl}
 \mathsf{Pr}_{\mathcal{S}_{sub}} & = & \left ( \mathsf{Id}_{|\mathcal{S}_{sub}}
\right)^\sharp & = &  \mathsf{Pr}_{\mathcal{S}_{sub}} \circ  \mathsf{Pr}_{\mathcal{S}_{sub}}.
\end{array}
$$
\normalsize
\end{prop}

\begin{defn}[Restricted mapping] \label{def_restriction}Let $f:\mathcal{E} \rightarrow \mathcal{F}$ be a
mapping and $\mathcal{A} \subseteq \mathcal{E}$. We will denote
$f_{|\mathcal{A}} : \mathcal{A} \rightarrow \mathcal{F} $ the
mapping defined by $f_{|\mathcal{A}} = f \circ \id_{|\mathcal{A}}$
where $\id_{|\mathcal{A}} : \mathcal{A} \rightarrow \mathcal{E}$
is the canonical injection from $\mathcal{A}$ to $\mathcal{E}$.
Similarly, let $\mathcal{B} \subseteq \mathcal{F}$ with $\im f
\subseteq \mathcal{B}$. Mapping $_{\mathcal{B}|}f : \mathcal{E}
\rightarrow \mathcal{B}$ is defined by $f = \id_{|\mathcal{B}}
\circ {_{\mathcal{B}|}f}$, where $\id_{|\mathcal{B}} : \mathcal{B}
\rightarrow \mathcal{F}$.
\end{defn}

\begin{prop} \label{prop:fcodomainedomaine}
Let $f : \mathcal{D} \rightarrow \mathcal{E}$ be a residuated
mapping and $\mathcal{D}_{sub}$ (resp. $\mathcal{E}_{sub}$) be a
complete subsemiring of $\mathcal{D}$ (resp. $\mathcal{E}$).
\begin{enumerate}
\item[$1.$]Mapping $f_{|\mathcal{D}_{sub}}$ is residuated and its residual is
given by : \small
$$ (f_{|\mathcal{D}_{sub}})^\sharp = ( f \circ \mathsf{Id}_{|\mathcal{D}_{sub}})^\sharp = \mathsf{Pr}_{S_{sub}} \circ
f^\sharp.
$$
\normalsize
\item[$2.$] If $\im f \subset \mathcal{E}_{sub}$ then mapping
$_{\mathcal{E}_{sub}|}f$ is residuated and its residual is given
by: \small
$$
\begin{array}{lcl}
\left ( _{\mathcal{E}_{sub}|} f \right )^\sharp & = & f^\sharp
\circ \mathsf{Id}_{|\mathcal{E}_{sub}} = \left ( f^\sharp \right
)_{|\mathcal{E}_{sub}}.
\end{array}
$$
\normalsize
\end{enumerate}

\end{prop}
\begin{pf}
Statement 1 follows directly from property \ref{prop:rescompo} and
proposition \ref{prop:canonicalinjection}. Statement 2 is obvious
since $f$ is residuated and $\im f \subset \mathcal{E}_{sub}
\subset \mathcal{E}$.\qed
\end{pf}

In \citep{cuning03}, the authors propose to compute a pair $(x,y)$
satisfying the following equation :
\begin{equation}\label{ax=by} a \otimes x= b \otimes
y.\end{equation} Theorems \ref{suite} and  \ref{thm:ax=by} given
below recall how to compute  such a pair, called simply a solution
hereafter.
\begin{thm}\label{suite}
Let $$ \Pi : x \mapsto L_a^\sharp \circ L_b \circ  L_b^\sharp
\circ L_a(x)= a \lslash (b \otimes (b \lslash (a \otimes x)))$$ be
a mapping defined over a semiring $\mathcal{S}$ and consider the following
iterative scheme :
\begin{itemize}
\item[] Let $x_0 \in \mathcal{S}$ be an arbitrary element
\item[] do $x_{n+1}=\Pi(x_n)$
\item[] until $x_{m+1}=x_{m}$ for $m \in \mathbb{N}$
\end{itemize}

If function $\Pi$ admits a fixed point $x$ such that $ x\preceq \Pi(x_0)$ and
$x\neq \varepsilon$ then the previous algorithm converges toward
the greatest finite fixed point in the following  principal order ideal
 $_\downarrow\mathcal{X}_{\Pi(x_0)}=\{x | x \preceq \Pi(x_0) \}$.
\end{thm}
\begin{pf}
First let us recall that a lower-bounded non increasing integer
sequence converges in a finite number of steps to the greatest
fixed point. Then if function $\Pi$ admits a fixed point $x$ such that $x\preceq \Pi(x_0)$
and $x\neq \varepsilon$ then it is sufficient
to show that the sequence $x_n$ ($n=1,2,..$) is non increasing.
Definition \ref{def_residuation} yields that $L_b \circ L_b^\sharp
\preceq \id$ and $L_a \circ L_a^\sharp \preceq \id$, then
$\Pi(x)\preceq L_a^\sharp \circ L_a(x)$ and
\begin{eqnarray}
\Pi \circ \Pi & = & L_a^\sharp \circ L_b \circ L_b^\sharp
\circ L_a \circ L_a^\sharp \circ L_b \circ  L_b^\sharp \circ L_a\\
\Pi \circ \Pi& \preceq & L_a^\sharp \circ L_b \circ L_b^\sharp
\circ L_a=\Pi \end{eqnarray}
 therefore $ x_{n+1} \preceq  x_n$, $i.e.$ the sequence is non increasing.\\
Obviously if $x_{n+1}=\Pi(x_n)$ then $x_{n+1}\preceq \Pi(x_0)$,
and it is the greatest in  $_\downarrow\mathcal{X}_{\Pi(x_0)}$.\qed
\end{pf}

\begin{thm}\label{thm:ax=by}
Define $y=L_b^\sharp \circ L_a(x)=b \lslash (a \otimes x)$. If the function $\Pi$ defined in theorem \ref{suite}
 admits a fixed point $x\neq \varepsilon$, then
the pair $(x,y)$  is a solution of (\ref{ax=by}).
\end{thm}
\begin{pf}
If $x=\Pi(x)$ and $y=L_b^\sharp \circ L_a(x)$ then, by using twice
equation (\ref{residuated1}), we obtain
$$
L_a(x)=L_a \circ  L_a^\sharp \circ L_b \circ  L_b^\sharp \circ
L_a(x) \preceq L_b \circ  L_b^\sharp \circ L_a(x)=L_b(y) \preceq
L_a(x).
$$
Therefore all these terms are equal, and $a \otimes x =
L_a(x)=L_b(y)= b \otimes y$. \qed
\end{pf}

\section{Semiring of intervals}

A semiring of interval may be constructed by considering a
semiring of pairs. The set of pairs $(x',x'')$ with $x' \in
\mathcal{S}$ and $x'' \in \mathcal{S}$ endowed with the following
coordinate-wise algebraic operations : \small
$$
\begin{array}{lcl}
& (x',x'') \oplus (y',y'') \triangleq (x' \oplus y' , x'' \oplus
y'') \\
 \textnormal{~~and~~} & (x',x'') \otimes (y',y'')
\triangleq (x' \otimes y' , x'' \otimes y''),
\end{array}
$$
\normalsize is a semiring denoted by $\cd$ with
$(\varepsilon,\varepsilon)$ as the zero element and $(e,e)$ as the
identity element (see definition \ref{def:semiring}).
Some of the preliminary statements of this
section are adapted from \citep{litvinov01c}.

\begin{rem}
\label{rem:order_rel}The operation $\oplus$ generates the
corresponding canonical partial order $\preceq_{\mathcal{C}}$ in
$\cd: \\ (x',x'') \oplus (y',y'') = (y',y'') \Leftrightarrow
(x',x'') \preceq_\mathcal{C} (y',y'') \Leftrightarrow x'
\preceq_\mathcal{S} y'$ and $x'' \preceq_\mathcal{S} y''$ where $
\preceq_\mathcal{S}$ is the order relation in $\mathcal{S}$.
\end{rem}

\begin{prop}
If the semiring $\mathcal{S}$ is complete, then the semiring $\cd$
is complete and its top element is given by $(\top,\top)$.
\end{prop}

\begin{notation} Consider the following mappings over
$\cd$ :
\small
$$
\begin{array}{lclr}
L_{(a',a'')} & : & (x',x'') \mapsto (a',a'') \otimes (x',x'') & ~~~~~~~~~~~~~~\textnormal{(left multiplication by $(a',a'')$)}\\
R_{(a',a'')} & : & (x',x'') \mapsto (x',x'') \otimes (a',a'') &
~~~~~~~~~~~~\textnormal{(right multiplication by $(a',a'')$)}.
\end{array}
$$
\normalsize
\end{notation}
\begin{prop} \label{res_sur C(D)}
The mappings $L_{(a',a'')}$ and $R_{(a',a'')}$ defined over $\cd$
are both residuated. Their residuals are equal to
$L_{(a',a'')}^\sharp (b',b'') = (a',a'') \lslash (b',b'')= (a'
\lslash b' , a'' \lslash b'')$ and $R_{(a',a'')}^\sharp (b',b'') =
(b',b'') \rslash (a',a'') = (b' \rslash a', b'' \rslash a'')$.
\end{prop}
\begin{pf}
Observe that $L_{(a',a'')} \left (\bigoplus_{(x',x'') \in X}
(x',x'') \right ) = \bigoplus_{(x',x'') \in X} L_{(a',a'')}\left (
x',x'' \right )$, (for every subset $X$ of
$\mathcal{C}(\mathcal{S}$)), moreover
$L_{(a',a'')}(\varepsilon,\varepsilon) = (a' \otimes \varepsilon,
a'' \otimes \varepsilon) = (\varepsilon,\varepsilon)$. Then
$L_{(a',a'')}$ is residuated (due to Theorem
\ref{thm:residuation}). Therefore, we have to find, for given
$(b',b'')$ and $(a',a'')$, the greatest solution $(x',x'')$ for
inequality $(a',a'') \otimes (x',x'') \preceq_\mathcal{C} (b',b'')
\Leftrightarrow (a' \otimes x' , a'' \otimes x'')
\preceq_\mathcal{C} (b', b'')$, moreover according to Remark
\ref{rem:order_rel} on the order relation induced by $\oplus$ on
$\cd$ we have, \small
$$
\begin{array}{lclclcl}
a' \otimes x' & \preceq_\mathcal{S} & b' & \textnormal{and} & a''
\otimes x'' & \preceq_\mathcal{S} & b''.
\end{array}
$$
\normalsize Since the mappings $x' \mapsto a' \otimes x'$ and $x''
\mapsto a'' \otimes x''$ are residuated over $\mathcal{S}$ (cf.
Corollary \ref{cor:eq_res}), we have $x' \preceq_\mathcal{S} a'
\lslash b'$ and $x'' \preceq_\mathcal{S} a'' \lslash b''$. Then,
we obtain $L_{(a',a'')}^\sharp(b',b'') = (a' \lslash b', a''
\lslash b'')$. \qed
\end{pf}
\begin{notation}
The set of pairs $(\widetilde{x}', \widetilde{x}'')$ s.t. $
\widetilde{x}'\preceq  \widetilde{x}''$ is denoted by
$\mathcal{C}_\mathtt{O}(\mathcal{S})$.
\end{notation}
\begin{prop} \label{prop:subsemiringcod}
Let $\mathcal{S}$ be a complete semiring. The set
$\mathcal{C}_\mathtt{O}(\mathcal{S})$ is a complete subsemiring of
$\cd$.
\end{prop}
\begin{pf}
Clearly $\mathcal{C}_\mathtt{O}(\mathcal{S}) \subset \cd$ and it
is closed for $\oplus$ and $\otimes$ since $\widetilde{x}' \oplus
\widetilde{y}' \preceq \widetilde{x}'' \oplus \widetilde{y}''$ and
$\widetilde{x}' \otimes \widetilde{y}' \preceq \widetilde{x}''
\otimes \widetilde{y}''$ whenever $\widetilde{x}' \preceq
\widetilde{x}''$ and $\widetilde{y}' \preceq \widetilde{y}''$.
Moreover zero element $(\varepsilon,\varepsilon)$, unit element
$(e,e)$ and top element $(\top,\top)$ of $\cd$ are in
$\mathcal{C}_\mathtt{O}(\mathcal{S})$. Therefore definition \ref{def:subsemiring} yields the result.\qed
\end{pf}
\begin{prop}
The canonical injection
$\mathsf{Id}_{|\mathcal{C}_\mathtt{O}(\mathcal{S})} :
\mathcal{C}_\mathtt{O}(\mathcal{S}) \rightarrow \cd$ is
residuated. Its residual
$(\mathsf{Id}_{|\mathcal{C}_\mathtt{O}(\mathcal{S})})^\sharp$ is a
projector denoted by
$\mathsf{Pr}_{\mathcal{C}_\mathtt{O}(\mathcal{S})} $. Its practical
computation is given by : \small
\begin{equation}
\begin{array}{lcl}
\mathsf{Pr}_{\mathcal{C}_\mathtt{O}(\mathcal{S})}((x',x'')) & = &
(x' \wedge x'', x'') = (\widetilde{x}', \widetilde{x}'').
\end{array}
\end{equation}
\normalsize
\end{prop}
\begin{pf}
It is a direct application of proposition
\ref{prop:canonicalinjection}, since
$\mathcal{C}_\mathtt{O}(\mathcal{S})$ is a subsemiring of $\cd$.
Practically, let $(x',x'') \in \cd$, we have
$\mathsf{Pr}_{\mathcal{C}_\mathtt{O}(\mathcal{S})}((x',x'')) =
(\widetilde{x}', \widetilde{x}'') = (x' \wedge x'', x'')$, which
is the greatest pair such that :

$$   \widetilde{x}' \preceq x', ~~~\widetilde{x}'' \preceq x'' ~~\textnormal{and}~~\widetilde{x}' \preceq \widetilde{x}''.$$\qed
\end{pf}

\begin{prop} \label{prop:residuation_couple_ordonnes}
Let $( \widetilde{a}', \widetilde{a}'') \in \cod$, then the mapping
${}_{{}_{{}_{\cod|}}}L_{{( \widetilde{a}',
\widetilde{a}'')}_{|\cod}} : \cod \rightarrow \cod$ is residuated.
Its residual is given by \small
$$
\begin{array}{lcl}
\left ( {}_{{}_{{}_{\cod|}}} L_{{( \widetilde{a}',
\widetilde{a}'')}_{|\cod}} \right )^\sharp & = &
\mathsf{Pr}_{\mathcal{C}_\mathtt{O}(\mathcal{S})} \circ \left (
L_{( \widetilde{a}', \widetilde{a}'')} \right )^\sharp \circ
\mathsf{Id}_{|\cod}.
\end{array}
$$
\normalsize
\end{prop}

\begin{pf}
Since  $(\widetilde{a}', \widetilde{a}'') \in \cod \subset \cd$,
it follows directly from proposition \ref{res_sur C(D)} that
mapping $L_{( \widetilde{a}', \widetilde{a}'')}$ defined over
$\mathcal{C}(\mathcal{S})$ is residuated. Furthermore, $\cod$
being closed for $\otimes$ we have $\im {L_{( \widetilde{a}',
\widetilde{a}'')|}}_{\cod} \subset \cod $, it follows from
definition \ref{def_restriction} and proposition
\ref{prop:fcodomainedomaine} that :
$$
\begin{array}{lclr}
\left ( {}_{{}_{{}_{\cod|}}}L_{{( \widetilde{a}',
\widetilde{a}'')}_{|\cod}} \right )^\sharp & = & \left ( L_{(
\widetilde{a}', \widetilde{a}'')} \circ \mathsf{Id}_{|\cod} \right
)^\sharp \circ \mathsf{Id}_{|\cod} &
 \\
& = & \mathsf{Pr}_{\mathcal{C}_\mathtt{O}(\mathcal{S})} \circ
\left ( L_{( \widetilde{a}', \widetilde{a}'')} \right )^\sharp
\circ \mathsf{Id}_{|\cod}. &
\end{array}
$$
\normalsize

Then, by considering $(\widetilde{b}',\widetilde{b}'') \in \cod
\subset \cd$, the greatest solution in $\cod$ of
$L_{(\widetilde{a}',\widetilde{a}'')}((\widetilde{x}',\widetilde{x}''))
= (\widetilde{a}',\widetilde{a}'') \otimes
(\widetilde{x}',\widetilde{x}'') \preceq
(\widetilde{b}',\widetilde{b}'')$ is
$L_{(\widetilde{a}',\widetilde{a}'')}^\sharp
((\widetilde{b}',\widetilde{b}'')) =
(\widetilde{x}',\widetilde{x}'') =
(\widetilde{a}',\widetilde{a}'') \lslash
(\widetilde{b}',\widetilde{b}'') =
\mathsf{Pr}_{\mathcal{C}_\mathtt{O}(\mathcal{S})}((\widetilde{a}'
\lslash \widetilde{b}' , \widetilde{a}'' \lslash \widetilde{b}''))
=(\widetilde{a}' \lslash \widetilde{b}' \wedge \widetilde{a}''
\lslash \widetilde{b}'', \widetilde{a}'' \lslash
\widetilde{b}'')$. \qed

\end{pf}

\begin{defn} \label{def:closedinterval}
A (closed) interval in semiring $\mathcal{S}$ is a set of the form
$\mathbf{x} = [\underline{x},\overline{x}] = \{t \in \mathcal{S} |
\underline{x} \preceq t \preceq \overline{x} \}$, where
$(\underline{x},\overline{x}) \in \cod$, $\underline{x}$
(respectively, $\overline{x}$) is said to be the lower
(respectively, upper) bound of the interval $\mathbf{x}$.
\end{defn}
\begin{prop}
The set of intervals, denoted by $\mathrm{I}(\mathcal{S})$,
endowed with the following coordinate-wise algebraic operations :
\small
\begin{equation}\label{eq:SumProductInterval}
\begin{array}{lcl}
\mathbf{x} \stackrel{-}{\oplus} \mathbf{y} \triangleq  \left
[\underline{x} \oplus \underline{y}, \overline{x} \oplus
\overline{y} \right ] & \textnormal{~~~and~~~} & \mathbf{x}
\stackrel{-}{\otimes} \mathbf{y} \triangleq \left [\underline{x}
\otimes \underline{y}, \overline{x} \otimes \overline{y} \right ]
\end{array}
\end{equation}
\normalsize is a semiring, where the interval
$\pmb{\varepsilon}=[\varepsilon,\varepsilon]$ (respectively,
$\mathbf{e}=[e,e]$) is zero (respectively, unit) element of
$\mathrm{I}(\mathcal{S})$. Moreover, the semiring
$\mathrm{I}(\mathcal{S})$ is isomorphic to $\cod$.
\end{prop}
\begin{pf}
First,  $\underline{x} \oplus \underline{y} \preceq \overline{x}
\oplus \overline{y}$ and $\underline{x} \otimes \underline{y}
\preceq \overline{x} \otimes \overline{y}$ whenever $\underline{x}
\preceq \overline{x}$ and $\underline{y} \preceq \overline{y}$,
then $\mathrm{I}(\mathcal{S})$ is closed with respect to the
operations $\stackrel{-}{\oplus},\stackrel{-}{\otimes}$. From
definition \ref{def:semiring}, it follows directly that it is a
semiring. Obviously, it is isomorphic to $\cod$, indeed let  $
\Psi : \cod \rightarrow \mathrm{I}(\mathcal{S}),
(\widetilde{x}',\widetilde{x}'') \mapsto
[\underline{x},\overline{x}]=[\widetilde{x}',\widetilde{x}'']$ be
the mapping which maps an interval to an ordered pair. Obviously
$\Psi^{-1}$ is well defined and the both mappings  are
homomorphisms.\qed
\end{pf}
\begin{rem}
Operations (\ref{eq:SumProductInterval}) give the tightest interval containing all results of the same operations to arbitrary elements of its interval operands.
\end{rem}
\begin{rem} \label{def:sommeinfinie}
 Let $\mathcal{S}$ be a complete semiring and $\{\mathbf{x}_{\alpha}\}$ be an infinite subset of
$\mathrm{I}(\mathcal{S})$, the infinite sum of elements of this
subset is : \small
$$ \overline{\bigoplus_{\alpha}} \mathbf{x}_{\alpha} = \left [ \bigoplus_{\alpha} \underline{x}_{\alpha}, \bigoplus_{\alpha} \overline{x}_{\alpha} \right
].
$$
\normalsize
\end{rem}
\begin{rem}
If $\mathcal{S}$ is a complete semiring then
$\mathrm{I}(\mathcal{S})$ is a complete semiring. Its top element
is given by
$\pmb{\top}=[\top,\top]$.\\
Note that if $\mathbf{x}$ and $\mathbf{y}$ are intervals in
$\mathrm{I}(\mathcal{S})$, then $\mathbf{x} \subset \mathbf{y}$
iff $\underline{y} \preceq \underline{x} \preceq \overline{x}
\preceq \overline{y}$. In particular, $\mathbf{x} = \mathbf{y}$
iff $\underline{x} =
\underline{y}$ and $\overline{x} = \overline{y}$. \\
An interval for which $\underline{x} = \overline{x}$ is called
\textit{degenerate}. Degenerate intervals allow to represent
numbers without uncertainty. In this case we identify $\mathbf{x}$
with its element by writing $\mathbf{x} \equiv x$.
\end{rem}

\begin{prop} \label{res_interval}
Mapping $L_{\mathbf{a}} : \mathrm{I}(\mathcal{S}) \rightarrow
\mathrm{I}(\mathcal{S}),\mathbf{x} \mapsto \mathbf{a}
\stackrel{-}{\otimes} \mathbf{x}$ is residuated. Its residual is
equal to $$L_{\mathbf{a}}^\sharp : \mathrm{I}(\mathcal{S})
\rightarrow \mathrm{I}(\mathcal{S}), (\mathbf{x}) \mapsto
\mathbf{a} \overline{\lslash} \mathbf{x}= [\underline{a}\lslash
\underline{x} \wedge \overline{a}\lslash\overline{x} ,
\overline{a} \lslash \overline{x}].$$ Therefore, $ \mathbf{a}
\overline{\lslash} \mathbf{b}$ is the greatest solution of
$\mathbf{a} \stackrel{-}{\otimes} \mathbf{x} \preceq \mathbf{b} $,
and the equality is achieved if $\mathbf{b} \in \im
L_{\mathbf{a}}$.
\end{prop}
\begin{pf}
Let $ \Psi : \cod \rightarrow \mathrm{I}(\mathcal{S}),
(\widetilde{x}',\widetilde{x}'') \mapsto
[\underline{x},\overline{x}]=[\widetilde{x}',\widetilde{x}'']$ be
the mapping which maps an interval to an ordered pair. This
mapping defines an isomorphism, since it is sufficient to handle
the bounds to handle an interval. Then the result follows directly
from proposition \ref{prop:residuation_couple_ordonnes}.
 \qed
\end{pf}
\begin{rem}
We would show in the same manner that mapping $R_\mathbf{a} :
\mathrm{I}(\mathcal{S}) \rightarrow \mathrm{I}(\mathcal{S}),
\mathbf{x} \mapsto \mathbf{x} \stackrel{-}{\otimes} \mathbf{a}$ is
residuated.
\end{rem}

\begin{rem}
These results show that it is not sufficient to consider
independently the  intervals bounds to compute the solution of the
equation $\mathbf{a}  \stackrel{-}\otimes \mathbf{x} \preceq
\mathbf{b}$. Therefore solution $\mathbf{(x,y)}$ of equality
$\mathbf{a} \stackrel{-}\otimes \mathbf{x}=\mathbf{b}
\stackrel{-}\otimes \mathbf{y}$ can not  be obtained by using the
algorithm proposed in theorem \ref{suite} independently for each
bounds. The following interval version needs to be used with the
following mapping
\begin{equation} \mathbf{\Pi : x \mapsto  a \stackrel{-}\lslash (b
\stackrel{-}\otimes (b \stackrel{-}\lslash (a \stackrel{-}\otimes
x)))}.
\end{equation}
 According to theorem \ref{thm:ax=by}, if $\mathbf{\Pi}$ admits a fixed point
 $\mathbf{x}\neq [
\varepsilon,\varepsilon] $, the pair $\mathbf{(x,y)}$ with
$\mathbf{y=b \stackrel{-} \lslash (a \stackrel{-}\otimes x)}$ is a
solution of equation $\mathbf{a} \stackrel{-}\otimes \mathbf{x}
=\mathbf{b} \stackrel{-}\otimes \mathbf{y}$ and obviously the
iterative scheme given in theorem \ref{suite} still valid.
\end{rem}

\begin{exmp}
Let $\mathbf{A}=\begin{pmatrix}[2,3] & [5,9]\\ [7,8] &
[3,6]\end{pmatrix}$ and $\mathbf{B}=\begin{pmatrix}[ 1,9] & [2,5]
& [3,4]\\ [1,13] & [3,10]& [9,10]
\end{pmatrix}$ be matrices with entries in $(max,+)$ semiring. By
considering interval version of algorithm proposed in theorem
\ref{suite}, with $\mathbf{x_0}=\begin{pmatrix}[4,7] \\
[3,5]\end{pmatrix}$ the convergence is achieved with the following
value of $\mathbf{x}=\begin{pmatrix}[4,7] \\ [2,2] \end{pmatrix}$,
and then by applying proposition \ref{res_interval} we obtain
$\mathbf{y}=\begin{pmatrix}[2,2 ] \\ [5,5] \\ [2,5]\end{pmatrix}$.
These computations may be obtained thanks to the $(max,+)$
toolboxes interfaced with Scilab (see \citep{hardouin06,lhommeau00a})
and this specific example may be obtained at \small
\url{http://www.istia-angers.fr/~hardouin/LAAA.html}.
\end{exmp}
\begin{exmp}
By considering these same matrices in $(max,min)$ algebra, the
same algorithm with the rules given in example
\ref{exmp:residuationmaxmin}, the vectors obtained are :
$x=\begin{pmatrix}[3,7] \\ [2,5] \end{pmatrix}$ and
$y=\begin{pmatrix}[5,5 ] \\ [7,7] \\ [7,7]\end{pmatrix}$.
\end{exmp}

\section{Conclusion}
This work shows that Residuation theory and resolution of matrix
systems of the form $\mathbf{A x \preceq B}$ over interval semirings is not so
straightforward as interval extension of simpler linear-algebraic
operations considered in \citep{litvinov01c}.
 Nevertheless, by considering the right
calculus rule the classical algorithm
(see \citep{cuning03}) remains efficient to solve some interval systems.
Example are given both in $(max,+)$ semiring and $(max,min)$
semiring but may be extended to other semirings such as semiring
of series introduced in \citep{BICOQ,gaubert92a} and can be useful to solve control
problems such as the one considered in \citep{lhommeau05,ouerghi06b}.

\bibliographystyle{plain}        

\end{document}